\documentclass[11pt, a4paper]{article}
\usepackage{amsmath}
\usepackage{amsthm}
\usepackage{amssymb}
\usepackage{array}
\usepackage{multirow}
\usepackage{graphicx}
\usepackage{epsfig}
\usepackage{mathrsfs}
\usepackage{hhline}
\usepackage{longtable}
\usepackage{array}
\usepackage{enumerate}
\usepackage{xcolor}
\usepackage{float}

\usepackage{graphicx}
\graphicspath{%
    {converted_graphics/}
    {/}
}

\newcolumntype{"}{@{\hskip\tabcolsep\vrule width 1pt\hskip\tabcolsep}}

\usepackage[mathlines]{lineno}
\modulolinenumbers[2]
\newcommand*\patchAmsMathEnvironmentForLineno[1]{%
\expandafter\let\csname old#1\expandafter\endcsname\csname #1\endcsname  \expandafter\let\csname oldend#1\expandafter\endcsname\csname end#1\endcsname  \renewenvironment{#1}%
{\linenomath\csname old#1\endcsname}%
{\csname oldend#1\endcsname\endlinenomath}}%
\newcommand*\patchBothAmsMathEnvironmentsForLineno[1]{%
\patchAmsMathEnvironmentForLineno{#1}%
\patchAmsMathEnvironmentForLineno{#1*}}%
\AtBeginDocument{%
\patchBothAmsMathEnvironmentsForLineno{equation}%
\patchBothAmsMathEnvironmentsForLineno{align}%
\patchBothAmsMathEnvironmentsForLineno{flalign}%
\patchBothAmsMathEnvironmentsForLineno{alignat}%
\patchBothAmsMathEnvironmentsForLineno{gather}%
\patchBothAmsMathEnvironmentsForLineno{multline}%
}

\def\Z{\mathbb{Z}}
\def\N{\mathbb{N}}

\definecolor{vividviolet}{rgb}{0.62, 0.0, 1.0}
\def\nt{\noindent}
\def\ms{\medskip}
\def\rsq{\hspace*{\fill}$\blacksquare$\medskip}

\newtheoremstyle{de}
  {10pt}          
  {10pt}  
  {\rm}  
  {}
  {\bf}  
  {. }    
  { }    
  {}     
\theoremstyle{de}

\newtheorem{example}{Example}[section]

\newtheorem{problem}{Problem}[section]

\newtheoremstyle{theorem}
  {10pt}          
  {10pt}  
  {\it}  
  {}
  {\bf}  
  {. }    
  { }    
  {}     
\theoremstyle{theorem}

\numberwithin{equation}{section}
\def\Z{\mathbb{Z}}
\def\N{\mathbb{N}}

\newtheorem{theorem}{Theorem}[section]
\newtheorem{lemma}[theorem]{Lemma}

\newtheorem{conjecture}{Conjecture}[section]
\numberwithin{equation}{section}

\begin{document}
\baselineskip18truept
\normalsize
\begin{center}
{\mathversion{bold}\Large \bf Complete characterization of s-bridge graphs with local antimagic chromatic number 2}

\bigskip
{\large  G.C. Lau{$^{a}$}, W.C. Shiu{$^{b}$}, R. Zhang{$^c$}, K. Premalatha{$^d$}, M. Nalliah{$^{e,}$}\footnote{Corresponding author.} }\\

\medskip

\emph{{$^a$}Faculty of Computer \& Mathematical Sciences,}\\
\emph{Universiti Teknologi MARA (Segamat Campus),}\\
\emph{85000, Johor, Malaysia.}\\
\emph{geeclau@yahoo.com}\\

\medskip

\emph{{$^b$}Department of Mathematics, The Chinese University of Hong Kong,}\\
\emph{Shatin, Hong Kong.}\\
\emph{wcshiu@associate.hkbu.edu.hk}\\

\medskip

\emph{{$^c$}School of Mathematics and Statistics, \\Qingdao University, Qingdao 266071 China.}\\
\emph{rx.zhang87@qdu.edu.cn} \\

\medskip

\emph{{$^d$}National Centre for Advanced Research in Discrete Mathematics,\\Kalasalingam Academy of Research and Education, Krishnankoil, India.}\\
\emph{premalatha.sep26@gmail.com}\\

\medskip

\emph{{$^e$}Department of Mathematics, School of Advanced Sciences, \\Vellore Institute of Technology, Vellore-632 014, India.} \\
\emph{nalliahklu@gmail.com}
\end{center}


\medskip
\begin{abstract}
An edge labeling of a connected graph $G = (V, E)$ is said to be local antimagic if it is a bijection $f:E \to\{1,\ldots ,|E|\}$ such that for any pair of adjacent vertices $x$ and $y$, $f^+(x)\not= f^+(y)$, where the induced vertex label $f^+(x)= \sum f(e)$, with $e$ ranging over all the edges incident to $x$.  The local antimagic chromatic number of $G$, denoted by $\chi_{la}(G)$, is the minimum number of distinct induced vertex labels over all local antimagic labelings of $G$. In this paper, we characterize $s$-bridge graphs with local antimagic chromatic number 2.

\medskip

\end{abstract}

\tolerance=10000
\baselineskip12truept
\def\qed{\hspace*{\fill}$\Box$\medskip}

\def\s{\,\,\,}
\def\ss{\smallskip}
\def\ms{\medskip}
\def\bs{\bigskip}
\def\c{\centerline}
\def\nt{\noindent}
\def\ul{\underline}
\def\lc{\lceil}
\def\rc{\rceil}
\def\lf{\lfloor}
\def\rf{\rfloor}
\def\a{\alpha}
\def\b{\beta}
\def\n{\nu}
\def\o{\omega}
\def\ov{\over}
\def\m{\mu}
\def\t{\tau}
\def\th{\theta}
\def\k{\kappa}
\def\l{\lambda}
\def\L{\Lambda}
\def\g{\gamma}
\def\d{\delta}
\def\D{\Delta}
\def\e{\epsilon}
\def\lg{\langle}
\def\rg{\tongle}
\def\p{\prime}
\def\sg{\sigma}
\def\to{\rightarrow}

\newcommand{\K}{K\lower0.2cm\hbox{4}\ }
\newcommand{\cl}{\centerline}
\newcommand{\om}{\omega}
\newcommand{\ben}{\begin{enumerate}}

\newcommand{\een}{\end{enumerate}}
\newcommand{\bit}{\begin{itemize}}
\newcommand{\eit}{\end{itemize}}
\newcommand{\bea}{\begin{eqnarray*}}
\newcommand{\eea}{\end{eqnarray*}}
\newcommand{\bear}{\begin{eqnarray}}
\newcommand{\eear}{\end{eqnarray}}

\section{Introduction}

\nt A connected graph $G = (V, E)$ is said to be {\it local antimagic} if it admits a {\it local antimagic edge labeling}, i.e., a bijection $f : E \to \{1,\dots ,|E|\}$ such that the induced vertex labeling $f^+ : V \to \Z$ given by $f^+(u) = \sum f(e)$ (with $e$ ranging over all the edges incident to $u$) has the property that any two adjacent vertices have distinct induced vertex labels. Thus, $f^+$ is a coloring of $G$. Clearly, the order of $G$ must be at least 3.  The vertex label $f^+(u)$ is called the {\it induced color} of $u$ under $f$ (the {\it color} of $u$, for short, if no ambiguous occurs). The number of distinct induced colors under $f$ is denoted by $c(f)$, and is called the {\it color number} of $f$. The {\it local antimagic chromatic number} of $G$, denoted by $\chi_{la}(G)$, is $\min\{c(f) \;|\; f\mbox{ is a local antimagic labeling of } G\}$. Clearly, $2\le \chi_{la}(G)\le |V(G)|$.  Throughout this paper, we shall use $a^{[n]}$ to denote a sequence of length $n$ in which all terms are $a$, where $n\ge 2$. For integers $1\le a < b$, we let $[a,b]$ denote the set of integers from $a$ to $b$.

\ms\nt A graph consisting of $s$ paths joining two vertices is called an {\it $s$-bridge graph}, which is denoted by $\th(a_1,\dots, a_s)$, where $s\ge 2$ and $1\le a_1\le a_2 \le \cdots \le a_s$ are the lengths of the $s$ paths.  For convenience, we shall let $\th_s = \th(a_1,a_2,\ldots,a_s)$ if there is no confusion. In this paper, we shall characterize $\th_s$ with $\chi_{la}(\theta_s)=2$.

\ms\nt The contrapositive of the following lemma in \cite[Lemma 2.1]{LSN-DMGT} or~\cite[Lemma 2.3]{LSN-IJMSI} gives a sufficient condition for a bipartite graph $G$ to have $\chi_{la}(G)\ge 3$.

\begin{lemma}[{\cite[Lemma 2.3]{LSN-IJMSI}}]\label{lem-2part} Let $G$ be a graph of size $q$. Suppose there is a local antimagic labeling of $G$ inducing a $2$-coloring of $G$ with colors $x$ and $y$, where $x<y$. Let $X$ and $Y$ be the sets of vertices colored $x$ and $y$, respectively. Then $G$ is a bipartite graph with bipartition $(X,Y)$ and $|X|>|Y|$. Moreover,
$x|X|=y|Y|= \frac{q(q+1)}{2}$.
\end{lemma}

\nt Clearly, $2\le \chi(\th(a_1,a_2,\ldots,a_s)) \le 3$ and the lower bound holds if and only if  $a_1\equiv \cdots \equiv a_s  \pmod{2}$. By Lemma~\ref{lem-2part}, we immediately have the following lemma.

\begin{lemma}\label{lem-th=2} For $s\ge 2$ and $1\le i\le s$, if $\chi_{la}((\th(a_1,a_2,\ldots,a_s)) =2$, then $a_i \equiv 0 \pmod{2}$. Otherwise, $\chi_{la}((\th(a_1,a_2,\ldots,a_s)) \ge 3$. \end{lemma}

\section{Main Result}

\nt In this section, we assume $\chi_{la}(\th_s)=2$. So by Lemma~\ref{lem-th=2}, $\th_s=\th(a_1, \dots, a_s)$ is bipartite and all $a_i$ are even. When $s=2$, $\th_s$ is a cycle, whose local antimagic chromatic number is 3. Thus $s\ge 3$.

\ms\nt Let $u$ and $v$ be the vertices of $\th_s$ of degree $s$.  We shall call the $2s$ edges incident to $u$ or else to $v$ as {\it end-edges}. An integer labeled to an end-edge is called an {\it end-edge label}. A path that starts at $u$ and ends at $v$ is called a {\it $(u,v)$-path.}

\ms\nt For integers $i$ and $d$ and positive integer $s$, let $A_s(i;d)$ be the arithmetic progression of length $s$ with common difference $d$ and first term $i$. We first have two useful lemmas.

\begin{lemma}\label{lem-AP}
Suppose $s, d\in\N$.
\begin{enumerate}[(a)]
\item For $i,j\in\Z$, the sum of the $k$-th term of $A_s(i;d)$ and that of $A_s(j;-d)$ is $i+j$ for $k\in[1,s]$; and the sum of the $k$-th term of $A_s(i;d)$ and the $(k-1)$-st term of $A_s(j;-d)$ is $i+j+d$ for $k\in[2,s]$.
\item If $0<|i_1-i_2|<d$, then $A_s(i_1; d)\cap A_s(i_2,\pm d)=\varnothing$.
\end{enumerate}
\end{lemma}
\begin{proof}  It is easy to obtain (a). We prove the contrapositive of (b). Suppose $A_s(i_1; d)\cap A_s(i_2,\pm d)\ne\varnothing$. Let $a\in A_s(i_1; d)\cap A_s(i_2,\pm d)$. Now,  $a=i_1+j_1d=i_2+j_2d$ for some integers $j_1, j_2$. Thus, $|i_1 - i_2| = d|j_2-j_1|\ge d$ if $j_2 \ne j_1$  or else $|i_1 - i_2|=0$ if $j_2= j_1$.
\end{proof}

\begin{lemma}\label{lem-AP-2} Suppose $\delta\in[0,n^2]\setminus\{2, n^2-2\}$ for some  integer $n\ge 2$. There is a subset $B$ of $A_n(1;2)$ such that the sum of integers in $B$ is $\delta$.
\end{lemma}
\begin{proof} If $\delta=0$, choose $B=\varnothing$. Suppose $1\le \delta
\le 2n-1$ and $\delta\ne 2$. If $\delta$ is odd, then choose $B=\{\delta\}$. If $\delta$ is even, then $\delta\ge 4$. We may choose $B=\{1,\delta-1\}$.

\nt Suppose $\delta > 2n-1$, then may choose a largest $k$ such that $\kappa=\sum\limits_{j=n-k+1}^{n} (2j-1)\le \delta$. Let $\tau=\delta-\kappa$. By the choice of $k$, $0\le \tau< 2n-2k-1$. There are 3 cases.

\begin{enumerate}[1.]
\item Suppose $\tau= 0$. $B=A_k(2n-2k+1;2)$ is the required subset.
\item Suppose $\tau$ is odd. $B=A_k(2n-2k+1;2)\cup\{\tau\}$ is the required subset.
\item Suppose $\tau$ is even. If $\tau\ge 4$, then we may choose $B=A_k(2n-2k+1;2)\cup\{\tau-1, 1\}$. If $\tau=2$, then $2=\tau<2n-2k-1$. We have $k\le n-2$. If $k\le n-3$, then choose $B=A_{k-1}(2n-2k+3;2)\cup\{2n-2k-1, 3, 1\}$. If $k=n-2$, then $\kappa=n^2-4$ and hence $\delta=n^2-2$ which is not a case.
\end{enumerate}
\end{proof}

\nt Suppose $A_1$ and $A_2$ be two sequences of length $n$. We combine these two sequences as a sequence of length $2n$, denoted $A_1\diamond A_2$, whose $(2i-1)$-st term is the $i$-th term of $A_1$ and the $(2i)$-th term is the $i$-th term of $A_2$, $1\le i\le n$.

\begin{theorem}\label{thm-chilath=2}  For $s\ge 3$, $\chi_{la}(\th_s) = 2$ if and only if $\th_s=K_{2,s}$ with even $s\ge 4$ or  the size $m$ of $\th_s$ is greater than $2s+2$ and $\th_s$ is one of the following graphs:
\begin{enumerate}[1a.]
\item[1.] $\th(4l^{[3l+2]}, (4l+2)^{[l]})$, $l\ge 1$;
\item[2a.] $\th(2l-2, (4l-2)^{[3l-1]})$, $l\ge 2$;
\item[2b.] $\th(2, 4^{[3]}, 6)$; $\th(4, 8^{[5]}, 10^{[2]})$; $\th(6, 12^{[7]}, 14^{[3]})$;
\item[3a.] $\th(4l-2-2t, 2t, (4l-4)^{[l]}, (4l-2)^{[l-2]})$, $2\le l\le t\le \frac{5l-2}{4}$;
\item[3b.] $\th(4l-2-2t, 2t-2, (4l-4)^{[l-1]}, (4l-2)^{[l-1]})$, $2\le l\le t\le \frac{5l}{4}$;
\item[4.] $\th(2t, 4s-6-2t, 2s-4, (4s-6)^{[s-3]})$, $\frac{2s-3}{8}\le t\le \frac{6s-5}{8}$, $s\ge 4$.
\end{enumerate}  \end{theorem}

\begin{proof} Note that $K_{2,s} = \th(2^{[s]})$.  In~\cite[Theorems 2.11 and 2.12]{Arumugam}, the authors obtained $$\chi_{la}(K_{2,s}) = \begin{cases} 2 & \mbox{ if $s\ge 4$ is even}, \\ 3 & \mbox{ otherwise.} \end{cases}$$  We only consider $\th_s\ne K_{2,s}, s\ge 3$.   Suppose $\chi_{la}(\th_s)=2$.  Since each $a_i$ is even, $\th_s$ has even size $m = \sum^s_{i=1} a_i \ge 2s + 2 \ge 8$ edges and order $m-s+2$. Let $f$ be a local antimagic labeling that induces a 2-coloring of $\th_s$ with colors $x$ and $y$. Without lost of generality, we may assume $f^+(u) = f^+(v) = y$. Let $X$ and $Y$ be the sets of vertices with colors $x$ and $y$, respectively. It is easy to get that $|Y| = m/2-s+2$ and $|X| = m/2$. By Lemma~\ref{lem-2part}, we have $x|X| = y|Y| = m(m+1)/2$. Hence, $x=m+1\ge 2s+3\ge 9$ is odd, $y = m(m+1)/(m-2s+4)$ and $y \ge (1 + 2 + \cdots + 2s)/2 = (2s^2+s)/2$.

\ms\nt Note that $\th_s$ has at least 2 adjacent non-end-edges. Suppose $z_1z_2$ is not an end-edge with $f(z_1z_2) = l$. Without loss of generality, we assume $f^+(z_1) = x$, $f^+(z_2) = y$. Since $z_1z_2$ is not an end-edge, there is another vertex $z_3$ such that $z_1z_2z_3$ forms a path. So, $f(z_2z_3) = y-l$. Since $1\le y-l\le m$, we have $l\ge y-m = y-x+1$. Consequently, all integers in $[1,y-x]$ must be assigned to end-edges. So, $y-x\le 2s$. Moreover, since $l\ne y-l$, we get $l\ne y/2$ so that $y/2$ must be an end-edge label when $y$ is even.

\ms\nt Solving for $m$, we get $m = \frac{1}{2}(y - 1 \pm \sqrt{y^2+14y-8ys+1})$. Hence, $y^2 + 14y - 8ys + 1 = t^2 \ge 0$, where $t$ is a nonnegative integer. This gives $(y+7-4s)^2 + 1 -(7-4s)^2 = t^2$ or $(y+7-4s-t)(y+7-4s+t)=8(s-2)(2s-3)$. By letting $a = y+7-4s-t$ and $b=y+7-4s+t$, we have $2y+14-8s=a+b$ with $ab=8(2s^2-7s+6)=8(s-2)(2s-3)$. Clearly, $b\ge a > 0$. Since $a,b$ must be of same parity, we have both $a,b$ are even.

\ms\nt Recall that $y - (2s^2+s)/2\ge 0$. Now
\begin{align}y -(2s^2+s)/2 & = 4s-7+\frac{a+b}{2}-\frac{2s^2+s}{2}\nonumber\\
& = \frac{a+b}{2} - \frac{2s^2-7s +6}{2}-4=\frac{a+b}{2} -\frac{ab}{16}-4\nonumber\\
& = \frac{8a+8b-ab-64}{16} = -\frac{(a-8)(b-8)}{16}. \label{eq-y}
\end{align}
This implies that $a\le 8$.

\ms\nt We shall need the following claim which is easy to obtain. Through out the proof, by symmetry, we always assume $\a_1 < \b_r$.

\ms\nt{\bf Claim:} {\it  Let $\phi$ be a labeling of a path $P_{2r+1}= v_1v_2\cdots v_{2r+1}$ with $\phi(v_{2i-1}v_{2i}) = \a_i$ and $\phi(v_{2i}v_{2i+1}) = \b_i$ for $1\le i\le r$. Suppose $\phi^+(v_{2j}) = x$ for $1\le j\le r$ and $\phi^+(v_{2k+1}) = y$ for $0\le k\le r$, where $y > x$, then $\a_1+\b_1=x$, $\{\a_1, \a_2, \ldots, \a_r\}$ is an increasing sequence with common difference $y-x$ while $\{\b_1, \b_2, \ldots, \b_r\}$ is a decreasing sequence with common difference $y-x$.  }

\ms\nt {\bf Case (1).} Suppose $a=8$. By \eqref{eq-y} we have $y = (2s^2+s)/2$  which implies $s$ is even. Express $t$ and $y$ in terms of $s$. This gives (i) $m=s^2-3s/2-1$ which implies $s\equiv 2\pmod{4}$ and $x=s^2-3s/2$ or (ii) $m=2s$. Since $m\ge 2s+2$, (ii) is not a case. In (i), $y-x = 2s$ so that all integers in $[1,2s]$ are end-edge labels.

\ms\nt Let $P$ be a $(u,v)$-path of $\th_s$ with length $2r$ whose end-edges are labeled by integers in $[1, 2s]$. Suppose one of its end-edges is labeled by $\a_1$. By the claim, another end-edge is labeled by $\b_r=\b_1-(r-1)(y-x)=x-\a_1-2rs+2s\le 2s$. So \[2r \ge \frac{x-\a_1}{s}\ge\frac{s^2-3s/2-2s}{s}=s-\frac{7}{2}.\]
Since $s$ and $2r$ are even,  $2r\ge s-2$. Since $\b_r\ge 2$, we have $2r\le \frac{1}{s}(x-\a_1+2s-2)< s+\frac{1}{2}$.
Thus, each $(u,v)$-path of $\th_s$ is of length $s$ or $s-2$. Suppose $\th_s$ has $h$ path(s) of length $s$ and $(s-h)$ path(s) of length $s-2$. We now have $sh + (s-h)(s-2) = m$. Therefore, $h = (s-2)/4$. Thus, $\th_s = \th((s-2)^{[(3s+2)/4]}, s^{[(s-2)/4]})$ for $s\equiv 2\pmod{4}$.

\ms\nt Let $s = 4l+2$, $l\ge 1$. We now show that $\th((s-2)^{[(3s+2)/4]}, s^{[(s-2)/4]})=\th((4l)^{[3l+2]}, (4l+2)^{[l]})$ admits a local antimagic 2-coloring. Recall that $m=16l^2+10l$, $x=16l^2+10l+1$, $y=16l^2+18l+5$ and $y-x=8l+4$.
\begin{enumerate}[{Step }1:]
\item Label the edges of the path $R_i$ of length $4l+2$ by using the sequence $A_{2l+1}(i; 8l+4)\diamond A_{2l+1}(x-i; -8l-4)$ in order, $1\le i\le l$.  Note that, as a set $A_{2l+1}(x-i;-8l-4) = A_{2l+1}(2l+1-i;8l+4)$. So, by Lemma~\ref{lem-AP}(b), $A_{2l+1}(i;8l+4)\diamond A_{2l+1}(x-i;-8l-4)$ for all $i\in [1, l]=U_1$ form a partition of $\bigcup\limits_{j=0}^{2l} [(8l+4)j+1, (8l+4)j+2l]$.  By Lemma~\ref{lem-AP}(a), we see that all induced labels of internal vertices are $x$ and $y$ alternatively. Now, integers in $[1, 2l]$ are end-edge labels.
\item Label  the edges of the path $Q_j$ of length $4l$ by the sequence $A_{2l}(\a;8l+4)\diamond A_{2l}(x-\a;-8l-4)$, where $\a$ is the $j$-th integer of the sequence $U_2=[3l+1,4l+1] \cup[4l+3,5l+1]\cup\{5l+3,6l+3\}\cup[7l+5,8l+4]$ in order, $1\le j\le 3l+2$. Note again, $A_{2l}(\a;8l+4)\diamond A_{2l}(x-\a;-8l-4)$ for all $\a \in  U_2$ form a partition of $\bigcup\limits_{j=0}^{2l-1} [(8l+4)j+2l+1, (8l+4)j+8l+4]$. By Lemma~\ref{lem-AP}(a), we see that all induced labels of internal vertices are $x$ and $y$ alternatively. Now, integers in $[2l+1, 8l+4]$ are end-edge labels.
\item We now merge the end-vertices with end-edge labels in $U_1\cup U_2$ to get the vertex $u$. We then merge the other end-vertices with end-edge labels in $[1,8l+4]\setminus (U_1\cup U_2)$ to get the vertex $v$. Clearly, both $u$ and $v$ have induced vertex label $y$.
\end{enumerate}

\nt Note that $\left(\bigcup\limits_{j=0}^{2l} [(8l+4)j+1, (8l+4)j+2l]\right)\cup \left(\bigcup\limits_{j=0}^{2l-1} [(8l+4)j+2l+1, (8l+4)j+8l+4]\right)=[1, 16l^2+10l]$. So the labeling defined above is a local antimagic $2$-coloring for $\th((4l)^{[3l+2]}, (4l+2)^{[l]})$.

\ms\nt {\bf Case (2).} Suppose $a=6$. Now, $b=\frac{4}{3}(s-2)(2s-3)$. By \eqref{eq-y} we have $y = 2s(2s - 1)/3$ and hence $s\equiv 0,2\pmod 3$. Similar to Case~(1), since $m\ge 2s+2\ge 8$, we must have $m=(4s^2-8s)/3$ and $s\ge 5$. Now $y-x = 2s-1$. So integers in $[1,2s-1]\cup \{y/2=(2s^2-s)/3\}$ are end-edge labels.

\ms\nt Note that there are $s-1$ paths in $\th_s$ with both end-edges labeled with integers in $[1,2s-1]$. Suppose $P_{2r+1}$ is one of these $s-1$ paths.  Since $\a_1 < \b_r$, we have $\a_1\in[1,2s-2]$.  Now, $\b_r = (x-\a_1) - (r-1)(y-x) \le 2s-1=y-x$. Since $x=(4s^2-8s+3)/3$ and $y-x=2s-1$, we have that
\[(2s-6)(2s-1)/3+1= (4s^2-14s + 9)/3  \le x-\a_1 \le r(y-x) = r(2s-1)\]
Thus $r>(2s-6)/3\ge \frac{4}{3}$, i.e., $r\ge 2$. Hence $\b_{r-1}$ is labeled at a non-end-edge so that $\b_{r-1} = (x-\a_1) - (r-2)(y-x) \ge 2s$. Therefore,
\[(r-2)(2s-1) \le x-\a_1 - 2s \le  (4s^2-14s)/3 =(2s-6)(2s-1)/3-2<(2s-6)(2s-1)/3.\]
Consequently, $r-2 < (2s-6)/3=2s/3 -2$, i.e., $r<2s/3$. Combining the aboves, we have
$2s/3-2<r<2s/3$ so that $2s-6<3r<2s$.
This implies that $3r\in[2s-5, 2s-1]$. Since $s\not\equiv 1\pmod 3$ we have the following two cases.

\begin{enumerate}[a)]
\item Consider $s=3l$, $l\ge 2$. Since $3r\equiv0\pmod{3}$, we have $3r=2s-3$, i.e., $r=2l-1$. Thus, the $s$-th path must have length $m - (3l-1)(4l-2) = 2l-2$. Consequently, $\th_{3l} =  \th(2l-2, (4l-2)^{[3l-1]})$.

 We now show that $\th_{3l} =  \th(2l-2, (4l-2)^{[3l-1]})$ admits a local antimagic 2-coloring. For $l=2$, $\th_6 = \th(2,6^{[5]})$ with induced labels $y=44,x=33$ and the paths have vertex labels\\[1mm]
    \centerline{$\begin{array}{lll}22,11; & 1,32,12,21,23,10;  & 3,30,14,19,25,8;\\ 4,29,15,18,26,7;  & 5,28,16,17,27,6;  &  9,24,20,13,31,2.\end{array}$}

     All the left (respectively right) end vertices are merged to get the degree 6 vertex with induced label 44.

\nt For $l\ge 3$, we apply the following steps.
\begin{enumerate}[{Step }1:]
\item Label the edges of the path $R_i$ of length $4l-2$ by the sequence $A_{2l-1}(i;6l-1)\diamond A_{2l-1}((6l-1)(2l-1)-i;-6l+1)$ in order, $1\le i\le 3l-1$.
\item Label the path $Q$ of length $2l-2$ by the sequence $A_{l-1}(6l-1; 6l-1)\diamond A_{l-1}((6l-1)(l-2);-6l+1)$ in order. By Lemma~\ref{lem-AP}, one may check that all integers in $[1, 4l(3l-2)]$ are assigned after the step.
\item If we merge the end-vertices with end-edge labels in $[1, 3l-1]\cup\{y/2\}$ as $u$, then the induced label of $u$ is $\frac{1}{2}(9l^2-3l)+(6l^2-l)=\frac{1}{2}(21l^2-5l)$. Clearly it is less than $y=12l^2-2l$. The difference is $\delta=\frac{l}{2}(3l+1)$.
\item Consider the set of differences of two end-edge labels in $R_i$, $1\le i \le 3l-1$, which is $D=\{1,3,\ldots,6l-3\}=A_{3l-1}(1;2)$.
Clearly $3< \delta<(3l-1)^2-3$. By Lemma~\ref{lem-AP-2} we have a subset $B$ of $D$ such that the sum of numbers in $B$ is $\delta$.
\item Label all end-edges incident to $u$ by  $([1, 3l-1]\setminus \{\frac{6l-1-i}{2}\;|\; i\in B\})\cup \{\frac{6l-1+i}{2}\;|\; i\in B\}\cup\{6l^2-l\}$.
\end{enumerate}
\nt We have a local antimagic $2$-coloring for $\th_{3l} =  \th(2l-2, (4l-2)^{[3l-1]})$.

\item Consider $s=3l-1$, $l\ge 2$. Now, $3r=2s-4$ or $2s-1$ so that $r\in\{2l-2,2l-1\}$. Note that $r\ge 2$.

 Let the path with an end-edge label $y/2=(2s^2-s)/3$ be of length $2q$. Since $y/2\notin [1, 2s-1]$ and we assume $\a_1<\b_q$, this means $\b_q=(2s^2-s)/3=(3l-1)(2l-1)$.

If $q=1$, then $\a_1+\b_1=x$. This implies $\a_1+(3l-1)(2l-1)=(2l-1)(6l-5)$ and hence $\a_1=6l^2-11l+4$. Since $\a_1\le 2s-1=6l-3$, we get $6l^2-17l+7=(2l-1)(3l-7)\le 0$. The only solution is $l=2$ so that $s=5$. Note that $q=l-1$.

Suppose $q\ge 2$. Now $\a_q+\b_q=x$ and $\a_q=\a_1+(q-1)(y-x)$ implies that $\a_1=x-\b_q-(q-1)(2s-1)\le (2s-1)$. So $x-\b_q\le q(2s-1)$. In terms of $l$, we have $(2l-1)(6l-5)-(3l-1)(2l-1)\le q(6l-3)$. Thus $3l-4\le 3q$. This implies $q\ge l-1$.  Also note that $\b_1=\b_q+(2s-1)(q-1)\le m=\frac{1}{3}(4s^2-8s)$. In terms of $l$ we will obtain $(6l-3)q\le 6l^2-5l$. This implies $q\le l-\frac{2l}{6l-3}<l$.  Thus, $q\le l-1$. Combining the aboves, we have $q=l-1$, as in $q=1$ above. 

\ms\nt Now, suppose there are $k$ paths of length $4l-4$ and $3l-2-k$ paths of length $4l-2$. We then have $(2l-2)+k(4l-4)+(3l-2-k)(4l-2)=4(3l-1)(l-1)=m$. Solving this, we get $k=2l-1$.  Consequently, $\th_{3l-1} = \th(2l-2, (4l-4)^{[2l-1]},(4l-2)^{[l-1]})$ for $l\ge 2$.

\nt Recall that $y=12l^2-10l+2$, $x=12l^2-16l+5$, $y-x=6l-3$. Using the claim, we now have the followings.
\begin{itemize}
\item Consider the $l-1$ path(s) of length $4l-2$. We have $\a_1 = i < \b_{2l-1} = x-i - (y-x)(2l-2) =2l-1-i$. So $1\le i\le l-1$. Thus, numbers in $[1, l-1]$ must  serve as $\a_1$ for these $l-1$ path(s). Hence numbers in $[l, 2l-2]$ must serve as $\b_{2l-1}$ for these $l-1$ path(s). Thus, numbers in $[1, 2l-2]$ are assigned to these $l-1$ paths.
\item Consider the $2l-1$ paths of length $4l-4$. We have $2l-1\le \a_1 = i < \b_{2l-2} = x-i-(y-x)(2l-3)=8l-4-i$. So $2l-1\le i\le 4l-3$. Thus, numbers in $[2l-1, 4l-3]$ must  serve as $\a_1$ for these $2l-1$ path(s). Hence numbers in $[4l-1, 6l-3]$ must serve as $\b_{2l-2}$ for these $2l-1$ path(s). Thus, numbers in $[2l-1, 6l-3]\setminus\{4l-2\}$ are assigned to these $2l-1$ paths.

\item Consider the path of length $2l-2$.  This path must have $\a_1 = 4l-2$ and  $\b_{l-1} = 6l^2-5l+1 = y/2$.
\end{itemize}

Since $y/2$ is assigned to an end-edge incident to $w$, say, at the path of length $2l-2$, we have
\[\frac{1}{2}(25l^2-25l+6)=\sum\limits_{i=1}^{l-1} i +\sum\limits_{j=2l-1}^{4l-3}j +(6l^2-5l+1)\le f^+(w)=12l^2-10l+2.\]
We get $l=2,3,4$, which implies $s=5,8,11$, respectively.

\ms\nt For $s=5$, we get $\th_5=\th(2,4^{[3]},6)$ with induced vertex labels $y=30$, $x=21$. The labels of the paths are \[15,6; \quad  3,18,12,9;  \quad  4,17,13,8;   \quad  7,14,16,5;  \quad  1,20,10,11,19,2.\]
\nt For $s = 8$, we get $\th_8 = \th(4,8^{[5]},10^{[2]})$ with induced vertex labels $y =80$, $x = 65$. The labels of the paths are
\[\begin{array}{lll} 40,25,55,10; & \quad  5,60,20,45,35,30,50,15;  &\quad 6,59,21,44,36,29,51,14;\\
7,58,22,43,37,28,50,13; & \quad 8,57,23,42,38,27,49,12;  & \quad  11,48,26,39,41,24,56,9; \\
1,64,16,49,31,34,46,19,61,4; & \quad 2,63,17,48,32,33,47,18,60,3. &\end{array}\]
\nt For $s=11$, we get $\th_{11} = \th(6,12^{[7]},14^{[3]})$ with induced vertex labels $y=154$, $x=133$. The labels of the paths are
{\small \[\begin{array}{lll} 77,56,98,35,119,14; &  \quad  7,126,28,105,49,84,70,63,91,42,112,21;\\
8,125,29,104,50,83,71,62,92,41,113,20; & \quad   9,124,30,103,51,82,72,61,93,40,114,19;\\     10,123,31,102,52,81,73,60,94,39,115,18; &\quad     11,122,32,101,53,80,74,59,95,38,116,17;\\
12,121,33,100,54,79,75,58,96,37,117,16; & \quad    13,120,34,99,55,78,76,57,97,36,118,15;\\
1,132,22,111,43,90,64,69,85,48,106,27,127,6;  &\quad    2,131,23,110,44,89,65,68,86,47,107,26,128,5;\\       4,129,25,108,46,87,67,66,88,45,109,24,130,3.\end{array}\]
}
\end{enumerate}

\ms\nt {\bf Case (3).} Suppose $a=4$. In this case, $b=2(2s^2-7s+6)$ and $2y+14-8s = 4s^2-14s+16$. So $y=2s^2-3s+1$. Similar to the previous cases, $m=2s^2-5s+2$ only. Hence $s$ is even, $x = 2s^2-5s+3$ and $y-x = 2s-2$. So integers in $[1,2s-2]$ must be assigned to $2s-2$ end-edges. Let the remaining two end-edges are labeled by $\gamma_1$ and $\gamma_2$. We have
$4s^2-6s+2=2y=f^+(u)+f^+(v)=\sum\limits_{i=1}^{2s-2} i +\gamma_1+\gamma_2=(s-1)(2s-1)+\gamma_1+\gamma_2$. Thus, $\gamma_1+\gamma_2=2s^2-3s+1=y$.

\ms\nt Suppose $\gamma_1$ and $\gamma_2$ are labeled at the end-edges of the same path of length $2q$. Without loss of generality, $\a_1=\gamma_1$ and $\b_q=\gamma_2$ so that
$y=\a_1+\b_q=\a_1+(x-\a_1)-(q-1)(y-x)$. We have $q(y-x)=0$ which is impossible.
Therefore, $\gamma_1$ and $\gamma_2$ are labeled at different paths. Thus, there are $s-2$ paths whose end-edges are labeled by integers in $[1, 2s-2]$ and exactly two paths, say $Q_i$ with an end-edge label in $[1,2s-2]$ and another end-edge label $\gamma_i\ge 2s-1$, $i=1,2$.

\ms\nt Suppose $P_{2r+1}$ is a path with both end-edges labeled with integers in $[1,2s-2]$. By the assumption $1\le \a_1 < \b_r\le 2s-2$ and the claim, we have $\b_r = (x-\a_1) - (r-1)(y-x) \le 2s-2$. So
\[(2s-2)(s-3) = 2s^2-8s+6 < 2s^2-7s+5 \le  x-\a_1 \le r(y-x) = r(2s-2).\]
Thus $r \ge s-2\ge 2$. So $\b_{r-1}$ is labeled at a non-end-edge. Therefore, $\b_{r-1} = (x-\a_1)-(r-2)(y-x)\ge 2s-1$.
We have
\[(r-2)(2s-2)\le x-\a_1-2s+1 \le 2s^2 - 7s + 3 < 2s^2-6s+4 = (2s-2)(s-2).\]
So $r < s$. Thus $r \in \{s-2,s-1\}$.

\ms\nt Suppose $Q_i$ is of length $2r_i$ whose end-edges are labeled by $\a_{1,i}\in[1, 2s-2]$ and $\b_{r_i,i}=\gamma_i$.
So $\b_{r_i,i}=\gamma_i=x-\a_{1,i}-(r_i-1)(y-x)$.  Since $\gamma_1+\gamma_2=2s^2-3s+1$ is odd, $\gamma_2\ge \frac{1}{2}(2s^2-3s+2)$ and $\gamma_1\le \frac{1}{2}(2s^2-3s)$. Now
\begin{align*}(r_2-1)(2s-2)& =x-\a_{1,2}-\gamma_2\le 2s^2-5s+3-1-\frac{1}{2}(2s^2-3s+2)\\&=(2s^2-7s+2)/2=[(2s-2)(s-2)-s-2]/2 < (2s-2)(s-2)/2.\end{align*}
We have $2r_2-2<s-2$ and hence $2r_2\le s-2$.

\ms\nt Now
$y = \gamma_1+\gamma_2 = 2x - \a_{1,1}-\a_{1,2} - (r_1+r_2-2)(y-x)$ or $(r_1+r_2-1)(2s-2)=(r_1+r_2-1)(y-x)=x - \a_{1,1}-\a_{1,2}$. Since $\a_{1,1}, \a_{1,2}\in[1,2s-2]$,
\begin{align*}(s-1)(2s-2)& >(s-1)(2s-2)-s-2=2s^2-5s=x-3\ge (r_1+r_2-1)(2s-2)\\& \ge x-(4s-5)=2s^2-9s+8=(s-4)(2s-2)+s>(s-4)(2s-2).\end{align*}
So $s>r_1+r_2>s-3$ or $2r_1+2r_2\in\{2s-2, 2s-4\}$.   Thus $2r_1+s-2\ge 2r_1+2r_2\ge 2s-4$. So we have $2r_1\ge s-2\ge 2r_2$. Since $2r_1+2r_2\le 2s-2$ and $2r_2\ge 2$, $2r_1\le 2s-4$.

\ms\nt Without loss of generality, we may always assume that $\g_1$ is labeled at the end-edge of $Q_1$ incident to $u$. Since $s\ge 4$ and $f^+(u)=y$, $\g_2$ must be labeled at the end-edge of $Q_2$ incident to $v$.  Suppose there are $k$ paths of length $2s-4$ and $s-k-2$ paths of length $2s-2$. Therefore, $2(r_1+r_2) + k(2s-4) + (s-k-2)(2s-2) = 2s^2-5s+2$. So $2(r_1+r_2)=s-2+2k$. For convenience, we write $s=2l$ for $l\ge 2$.
\begin{enumerate}[(a)]
\item Suppose $2r_1+2r_2=4l-2$. Now, $k=l$ and $\th_{2l} = \th(4l-2-2r_1, 2r_1, (4l-4)^{[l]}, (4l-2)^{[l-2]})$ for $l\le r_1 \le 2l-2$. Since $l-1\ge r_2=2l-1-r_1$, $r_1\ge l$. Rewriting $r_1$ as $t$ we have $\th_{2l}=\th(4l-2-2t, 2t, (4l-4)^{[l]}, (4l-2)^{[l-2]})$ for $l\le t\le 2l-2$. Here $Q_2$ and $Q_1$ are $(u,v)$-paths of length $4l-2-2t$ and $2t$, respectively.

    Following we consider all $(u,v)$-paths of $\th_s$. Let the $(u,v)$-paths of length $4l-4$ be $R_i$, $1\le i\le l$ and the $(u,v)$-path(s) of length $4l-2$ be $T_j$, $1\le j\le l-2$. Let $T_{l-1}$ be the path obtained from $Q_2$ and $Q_1$ by merging the vertex $v$ of $Q_2$ and the vertex $u$ of $Q_1$. Hence $T_{l-1}$ is a $(u,v)$-path of length $4l-2$. Under the labeling $f$, the end-edge labels are in $[1, 4l-2]$ and the induced vertex labels of all internal vertices of $T_{l-1}$ are $x$ and $y$ alternatively.

\item Suppose $2r_1+2r_2=4l-4$. Now, $2r_1=4l-4-2r_2\le 4l-6$ so that $k=l-1$ and\\ $\th_{2l} = \th(4l-4-2r_1, 2r_1, (4l-4)^{[l-1]}, (4l-2)^{[l-1]})$ for $l-1 \le r_1 \le 2l-3$. Rewriting $r_1$ as $t-1$ we have
$\th_{2l} = \th(4l-2-2t, 2t-2, (4l-4)^{[l-1]},(4l-2)^{[l-1]})$ for $l \le t\le 2l-2$. Here $Q_2$ and $Q_1$ are $(u,v)$-paths of length $4l-2-2t$ and $2t-2$, respectively.

Following we consider all $(u,v)$-paths of $\th_s$. Let the path(s) of length $4l-4$ be $R_i$, $1\le i\le l-1$ and the path(s) of length $4l-2$ be $T_j$, $1\le j\le l-2$. Let $R_{l}$ be the path obtained from $Q_2$ and $Q_1$ by merging the vertex $v$ of $Q_2$ and the vertex $u$ of $Q_1$. Hence $R_{l}$ is a $(u,v)$-path of length $4l-4$. Under the labeling $f$, the end-edge labels are in $[1, 4l-2]$ and the induced vertex labels of all internal vertices of $R_{l}$ are $x$ and $y$ alternatively.
\end{enumerate}
For each case, after the merging, we have $l$ paths $R_i$ of length $4l-4$, $1\le i\le l$ and $l-1$ paths $T_j$ of length $4l-2$, $1\le j\le l-1$, where $l\ge 2$. All the end-edge labels are in $[1, 4l-2]$ under the labeling $f$.
Consider the $(u,v)$-path $R_i$ of length $2s-4=4l-4$. Suppose $x_i=\a_1$ is an end-edge label, then another end-edge label is $\b_{s-2}=(x-\a_1)-(s-3)(2s-2)\le 2s-2$. We have $\a_1\ge s-1$. By symmetry, $\b_{s-2}\ge s-1$. So all the $l$ paths $R_i$  have their end-edges labeled by integers in $[2l-1, 4l-2]$. Thus, all $(u,v)$-paths $T_j$  have their end-edges labeled by integers in $[1, 2l-2]$.

\ms\nt Let the label assigned to the end-edge of $T_j$ incident to $u$ be $y_j$.
\begin{enumerate}[(a)]
\item For the case $\th_{2l} = \th(4l-2-2t, 2t, (4l-4)^{[l]}, (4l-2)^{[l-2]})$, $2\le l \le t \le 2l-2$, $\g_1$ is the $(4l-2-2t+1)$-st edge label of $T_{l-1}$ so that $\g_1=y_{l-1}+(2l-1-t)(4l-2)$. Hence
    \[(4l-1)(2l-1)=f^+(u)=\g_1+ \sum_{j=1}^{l-1} y_j +\sum_{i=1}^{l} x_i=y_{l-1}+(2l-1-t)(4l-2)+\sum_{j=1}^{l-1} y_j +\sum_{i=1}^{l} x_i.\]
    We have
    \begin{align*}(2l-1-t)(4l-2)& =(4l-1)(2l-1)-y_{l-1}-\sum_{j=1}^{l-1} y_j -\sum_{i=1}^{l} x_i \\ & \ge (4l-1)(2l-1) -(2l-2)-\frac{(l-1)(3l-2)}{2}-\frac{l(7l-3)}{2}=3l^2-4l+2.
    \end{align*}
    This means\\ $t(4l-2)\le 2(2l-1)^2-(3l^2-4l+2)=5l^2-4l= \frac{1}{4}[(5l-1)(4l-2)-2l-2]<\frac{1}{4}(5l-1)(4l-2)$. Therefore, $t< \frac{5l-1}{4}$, i.e., $t\le \frac{5l-2}{4}$. Thus, $l\le t\le \frac{5l-2}{4}$.
\item For the case $\th_{2l}=\th(4l-2-2t, 2t-2, (4l-4)^{[l-1]},(4l-2)^{[l-1]})$ for $2\le l \le t\le 2l-2$, similarly we have
     \begin{align*}(2l-1-t)(4l-2)& =(4l-1)(2l-1)-x_{l}-\sum_{j=1}^{l-1} y_j -\sum_{i=1}^{l} x_i \\ & \ge (4l-1)(2l-1) -(4l-2)-\frac{(l-1)(3l-2)}{2}-\frac{l(7l-3)}{2}=3l^2-6l+2.
    \end{align*}
    This means\\ $t(4l-2)\le 2(2l-1)^2-(3l^2-6l+2)=5l^2-2l= \frac{1}{4}[(5l+1)(4l-2)-2l+2]<\frac{1}{4}(5l+1)(4l-2)$.  Therefore, $t< \frac{5l+1}{4}$, i.e., $t\le \frac{5l}{4}$. Thus, $l\le t\le \frac{5l}{4}$.
\end{enumerate}
Consequently, we have the following two cases.
\begin{enumerate}[(a)]
\item $\th_{2l} = \th(4l-2-2t, 2t, (4l-4)^{[l]}, (4l-2)^{[l-2]})$ for $2\le l \le t\le \frac{5l-2}{4}$, or else
\item $\th_{2l} = \th(4l-2-2t, 2t-2, (4l-4)^{[l-1]},(4l-2)^{[l-1]})$ for $2\le l \le t\le \frac{5l}{4}$.
\end{enumerate}

\ms\nt Now, we are going to find a local antimagic $2$-coloring for the above graphs.

\begin{enumerate}[(a)]
\item $\th_{2l} = \th(4l-2-2t, 2t, (4l-4)^{[l]}, (4l-2)^{[l-2]})$ for $2\le l \le t\le \frac{5l-2}{4}$.
\begin{enumerate}[Step 1:]
\item Label the edges of $T_j$ by the sequence $A_{2l-1}(l-1+j; 4l-2)\diamond A_{2l-1}(x-l+1-j; -4l+2)$, $1\le j\le l-1$. Note that we choose $\a_1=l-1+j$. This gives $\b_{2l-1}=l-j$. So, as a set $A_{2l-1}(x-(l-1+j); -4l+2)=A_{2l-1}(l-j; 4l-2)$. Thus, integers in $[1, 2l-2]$ are end-edge labels of all path(s) $T_j$ and integers in $\bigcup\limits_{j=1}^{l-1} [(j-1)(4l-2)+1, (j-1)(4l-2)+(2l-2)]$ are assigned.

\item Label the edges of the $(u,v)$-path $R_i$ by the sequence $A_{2l-2}(2l-2+i; 4l-2)\diamond A_{2l-2}(x-2l+2-i; -4l+2)$, $1\le i\le l$.
Note that we choose $\a_1=2l-2+i$. This gives $\b_{2l-2}=6l-3-(2l-2+i)=4l-1-i$. So, as a set $A_{2l-2}(x-2l+2-i; -4l+2)=A_{2l-2}(4l-1-i; 4l-2)$. Thus, integers in $[2l-1, 4l-2]$ are end-edge labels of all path(s) $R_i$ and integers in $\bigcup\limits_{i=1}^{l} [(i-1)(4l-2)+(2l-1), (i-1)(4l-2)+(4l-2)]$ are assigned. The set of difference between the two end-edge labels of a path $R_i$ is $D_2=\{1, 3, \dots, 2l-1\}=A_{l}(1;2)$.
\item Pick the $(u,v)$-path $T_{l-1}$ and separate it into two paths. Note that the end-edge labels of $T_{l-1}$ are $2l-2$ and $1$. The first $4l-2-2t$ edges form a $(u,v)$-path $Q_2$ and the remaining $2t$ edges form a $(u,v)$-path $Q_1$. Note that the label of $(4l-1-2t)$-th edge of $T_{l-1}$ is  $\g_1=(2l-1-t)(4l-2)+(2l-2)$.
\end{enumerate}
Thus, the above labeling is a local antimagic labeling. Under this labeling, the induced vertex label of $u$ is
\begin{align*}\sum_{j=1}^{l-1} (l-1+j) + \sum_{i=1}^{l} (2l-2+i) +\g_1 & = \frac{(l-1)(3l-2)}{2}+\frac{l(5l-3)}{2}+  (2l-1-t)(4l-2)+(2l-2)\\ &=  12l^2+2t-10l-4lt+1.\end{align*}
The difference from $y=8l^2-6l+1$ is $\delta(t)=4lt+4l-4l^2-2t=(4l-2)(t-l)+2l$. Clearly $2<\d(t)\le (4l-2)\frac{l-2}{4}+2l\le l^2$. Suppose $\d(t)=l^2-2$, then  $t=\frac{5l^2-4l-2}{4l-2}=\frac{5l-2}{4}+\frac{l-6}{2(4l-2)}$. Since $t\le \frac{5l-2}{4}$, $2\le l\le 6$. Since $t\in\Z$, $l=6$ and hence $t=7$. Thus, by Lemma~\ref{lem-AP-2}, we may choose $B\subset D_2$ to obtain a local antimagic $2$-coloring of $\th(4l-2-2t, 2t, (4l-4)^{[l]}, (4l-2)^{[l-2]})$ for $2\le l \le t\le \frac{5l-2}{4}$ and $(l,t)\ne (6,7)$. We shall provide a local antimagic $2$-coloring for the special case $(l,t)=(6,7)$ in Example~\ref{ex-3-1}(a)(ii).

\item $\th_{2l} = \th(4l-2-2t, 2t-2, (4l-4)^{[l-1]},(4l-2)^{[l-1]})$ for $2\le l \le t\le \frac{5l}{4}$.
\begin{enumerate}[Step 1:]
\item Label the edges of $T_j$ by the sequence $A_{2l-1}(j; 4l-2)\diamond A_{2l-1}(x-j; -4l+2)$, $1\le j\le l-1$.
   The set of difference between the last label and the first label of a paths $T_j$'s is $D_1=\{1, 3, \dots, 2l-3\}=A_{l-1}(1;2)$.
\item Label the edges of $R_i$ by the sequence $A_{2l-2}(3l-2+i; 4l-2)\diamond A_{2l-2}(x-3l+2-i; -4l+2)$, $1\le i\le l$. The set of difference between the last label and the first label of a paths $R_i$'s, $1\le i\le l-1$ is $D_2=\{-1, -3, \dots, -(2l-3)\}=A_{l-1}(-1;-2)$.
\item Pick the $(u,v)$-path $R_{l}$ and separate it into two paths. Note that the end-edge labels of $R_l$ are $4l-2$ and $2l-1$. The first  $4l-2-2t$ edges form a $(u,v)$-path $Q_2$ and the remaining  $2t-2$ edges form a $(u,v)$-path $Q_1$. Note that the label of $(4l-1-2t)$-th edge of $R_l$ is $\g_1=(2l-1-t)(4l-2)+(4l-2)$.
\end{enumerate}
Similar to the previous case, the above labeling is a local antimagic labeling. Under this labeling, the induced vertex label of $u$ is
\begin{align*} \sum_{j=1}^{l-1} j + \sum_{i=1}^{l} (3l-2+i) +\g_1 & = \frac{(l-1)l}{2}+\frac{l(7l-3)}{2}+(2l-1-t)(4l-2)+(4l-2)\\ &=12l^2+2t-6l-4lt.\end{align*}
The difference from $y=8l^2-6l+1$ is $\delta(t)=-4l^2-2t+4lt+1$. Clearly $\delta(t)$ is an increasing function of $t$. It is easy to show that $3\le 2l-1\le\delta(t)\le l^2-\frac{5l}{2}+1\le (l-1)^2-1$ when $l+ 1\le t\le \frac{5l}{4}$. We need to show that $\delta(t)\ne (l-1)^2-2$. Now $\delta((5l-1)/4))=\frac{2l^2-7l+3}{2}=(l-1)^2-\frac{3l-1}{2}<(l-1)^2-2$. If $\frac{5l}{4}\in\Z$, then $l\ge 4$. So $\delta(5l/4)=\frac{2l^2-5l+1}{2}=(l-2)^2-\frac{l+1}{2}<(l-1)^2-2$.
Thus $3\le\delta(t)\le l^2-\frac{5l}{2}+1\le (l-1)^2-2$ when $l+ 1\le t\le \frac{5l}{4}$. By Lemma~\ref{lem-AP-2}, we may choose $B\subset D_1$ and then we obtain a local antimagic $2$-coloring for $\th(4l-2-2t, 2t-2, (4l-4)^{[l-1]},(4l-2)^{[l-1]})$ for $l+1 \le t\le \frac{5l}{4}$.

The remaining case is $t=l$. For this case, $\delta(l)=-2l+1$. If $l\ne 3$, then we may choose $B=\{-(2l-3), -3, 1\}\subset D_1\cup D_2$.  When $l=3$, we have $t=3$. This is a special case with solution given in  Example~\ref{ex-3-2}(b).

\end{enumerate}

\ms\nt  {\bf Case (4).} Suppose $a=2$. In this case, $b=4(2s^2-7s+6)$ and $2y+14-8s=8s^2-28s+26$. So $y=4s^2-10s+6$. Similar to the previous cases we have $m=4s^2-12s+8$. Hence $x=4s^2-12s+9$.

\ms\nt Suppose $s=3$. We get $m=8$, $x = 9$ and $y = 12$. Thus, $\th_3=\th(2,2,4)$. The sequences we can use are $3,6$; $1,8$ and $4,5,7,2$ or else $3,6$; $1,8,4,5$ and $7,2$, both of which give no solution. We now assume $s\ge 4$.

\ms\nt Note that $y-x=2s-3$, $y$ is even and $y/2>2s-3$. Recall that if $y$ ie even, then $y/2$ is an end-edge label. Thus, integers in $[1,2s-3]\cup \{y/2\}$ are end-edge labels.

\ms\nt There are only 3 end-edge labels greater than $2s-3$. So there are at least $s-3$ paths with both end-edges labeled by integers in $[1,2s-3]$. Suppose $P_{2r+1}$ is one of these $s-3$ paths. Keep the notation defined in the claim and the assumption $\a_1 < \b_r$. So, $\a_1 \in [1,2s-4]$.

\ms\nt Now $\b_r = (x-\a_1) - (r-1)(y-x) \le 2s-3$. Since $x=4s^2-12s+9$ and $y-x=2s-3$, we have
\[ (2s-3)(2s-4) < 4s^2 - 14s + 13 \le  x-\a_1 \le r(y-x) = r(2s-3) \]
Thus, $r\ge 2s-3$.

\ms\nt Since $r\ge 4$, $\b_{r-1}$ is labeled at a non-end-edge. So $\b_{r-1} = (x-\a_1)-(r-2)(y-x)\ge 2s-2$ so that
\[ (r-2)(2s-3) \le x - \a_1 - 2s + 2 \le 4s^2 - 14s + 10 < (2s-3)(2s - 4). \]
So $r-2\le 2s-5$ or $r\le 2s-3$. Thus, $r=2s-3$. Note that, $\b_{2s-3}=2s-3-\a_1$.

\ms\nt Suppose $y/2=2s^2-5s+3$ is labeled at an end-edge of a path $Q$. Let the length of $Q$ be $2q$. So we have $\a_1 \le 2s-3$, $\b_q = y/2$  and $\b_1 = y/2 + (q-1)(2s-3)$. Now
$x=\a_1+\b_1=\a_1+y/2 +(q-1)(y-x)$ so that $2x> y+(2q-2)(y-x)$. We have $(2s-3)^2=x>(2q-1)(y-x)=(2q-1)(2s-3)$. Thus $2q-1<2s-3$, i.e., $q\le s-2$.

\ms\nt On the other hand, $2x=2\a_1+y+(2q-2)(y-x)\le 2(2s-3)+y+(2q-2)(y-x)=y+2q(y-x)$ so that $(2s-3)^2=x\le (2q+1)(y-x)=(2q+1)(2s-3)$. This means $2q+1\ge 2s-3$, i.e., $q\ge s-2$. Thus $q=s-2$. Consequently, $\th_s$ contains a path of length $2s-4$ with an end-edge label $\b_{s-2} = 2s^2-5s+3 = y/2$ so that $\a_i = i(2s-3)$ and $\b_i = 4s^2 - 14s + 12 - (i-1)(2s-3) = (2s-3)(2s-3-i)\ge (2s-3)(s-1)$ for $1\le i\le s-2$.

\ms\nt Let the remaining two end-edge labels be $\gamma_1$ and $\gamma_2$. Thus, $2y=f^+(u)+f^+(v)=\gamma_1+\gamma_2+y/2+(2s-3)(s-1)$. So $\gamma_1+\gamma_2=4s^2-10s+6=y$.

\nt Suppose $\gamma_1$ and $\gamma_2$ are labeled at the same path of length $2q$. By a similar proof of Case (3), we have
$4s^2-10s+6=\gamma_1+\gamma_2=\gamma_1+(x-\gamma_1)-(q-1)(y-x)=4s^2-12s+9-(q-1)(2s-3)$ which is impossible.

\ms\nt As a conclusion, there are exactly $s-3$ paths of length $4s-6$ whose end-edges are labeled by integers in $[1, 2s-4]$, one path of length $2s-4$ whose end-edges are labeled by $2s-3$ and $y/2$, two paths $Q_i$ of length $s_i$ whose end-edges are labeled by $\a_{1,i}\in [1, 2s-4]$ and $\gamma_i$, $i=1,2$. By counting the number of edges of the graph, we have $s_1+s_2=4s-6$. Thus, $\th_s=\th(2t,4s-6-2t,2s-4,(4s-6)^{[s-3]})$ for some $t\ge 1$.

Let us rename all $(u,v)$-paths.
\begin{itemize}
\item Let $R_1, \dots R_{s-3}$ be the $(u,v)$-paths in $\th_s$ of length $4s-6$. Let the end-edge label of $R_i$ incident to $u$ be $x_i$, $1\le i\le s-3$.

\item Let $P$ be the $(u,v)$-path of length $2s-4$ whose end-edge labels are $2s-3$ and $(s-1)(2s-3)$.

\item Let $Q_1$ be $(u,v)$-path of length $4s-6-2t$ whose end-edge labels are $\gamma_1$ and $x_{s-1}$. Let $Q_2$ be $(u,v)$-path of length $2t$ whose end-edge labels are $x_{s-2}$ and $\gamma_2$. Without loss of generality, we may assume that $\gamma_1<\gamma_2$. Since $\gamma_1+\gamma_2=y$, $\gamma_1<y/2<\gamma_2$. Also, without loss of generality, we may always assume that $\gamma_1$ is labeled at the end-edge incident to $u$. Thus, $x_{s-2}$ is labeled at the end-edge of $Q_2$ incident to $u$.

    Let $R_{s-2}$ be the labeled $(u,v)$-path obtained from $Q_2$ and $Q_1$ by merging the end vertex $v$ of $Q_2$ with the end vertex $u$ of $Q_1$. Therefore, $R_{s-2}$ satisfies the assumption of the Claim. Thus $x_{s-2}$ is labeled at the end-edge of $R_{s-2}$ incident to $u$. Now $\gamma_1=t(2s-3)+x_{s-2}$.
\end{itemize}

\nt Suppose $2s-3$ is labeled at the end-edge of $P$ incident to $u$, then
\begin{align*}& \quad\ 2(s-1)(2s-3) =f^+(u) =\sum\limits_{i=1}^{s-3} x_i +(2s-3)+x_{s-2}+\gamma_1\\
&=\sum\limits_{i=1}^{s-2} x_i +(2s-3)+[t(2s-3)+x_{s-2}]= \sum\limits_{i=1}^{s-2} x_i +(t+1)(2s-3)+x_{s-2}
\end{align*}
This means $(2s-t-3)(2s-3)=x_{s-2}+\sum\limits_{i=1}^{s-2} x_i\le (2s-4)+\frac{(s-2)(3s-5)}{2}$. Since $1\le t\le s-2$, $(s-1)(2s-3)\le (2s-4)+\frac{(s-2)(3s-5)}{2}=\frac{3s^2-7s+2}{2}$ which is impossible.
Thus, $(s-1)(2s-3)$ must be a label of the end-edge of $P$ incident to $u$. Consequently, we have
\begin{align*}&\quad \ 2(s-1)(2s-3)=f^+(u) =\sum\limits_{i=1}^{s-3} x_i +(s-1)(2s-3)+x_{s-2}+\gamma_1\\
&=\sum\limits_{i=1}^{s-2} x_i +(s-1)(2s-3)+[t(2s-3)+x_{s-2}]= \sum\limits_{i=1}^{s-2} x_i +(s-1+t)(2s-3)+x_{s-2}
\end{align*}
This means $(s-t-1)(2s-3)=x_{s-2}+\sum\limits_{i=1}^{s-2} x_i\ge 1+\frac{(s-2)(s-1)}{2}=\frac{s^2-3s+4}{2}=\frac{(2s-3)^2}{8}+\frac{7}{8}>\frac{(2s-3)^2}{8}$. Solve this inequality we have $t<\frac{6s-5}{8}$.

\nt Similarly, we have $(s-t-1)(2s-3)=x_{s-2}+\sum\limits_{i=1}^{s-2} x_i\le \frac{3s^2-7s+2}{2}=\frac{(6s-5)(2s-3)}{8}-\frac{7}{8}<
\frac{(6s-5)(2s-3)}{8}$. Solve this inequality we have $t>\frac{2s-3}{8}$.

\nt Hence
\[t\in \begin{cases} [2j-1, 6j-4] & \mbox{if }s=8j-4;\\
[2j-1, 6j-3] & \mbox{if } s=8j-3;\\
[2j, 6j-3] & \mbox{if } s=8j-2;\\
[2j, 6j-2] & \mbox{if }s=8j-1;\\
[2j, 6j-1] & \mbox{if }s=8j;\\
[2j, 6j] & \mbox{if } s=8j+1;\\
[2j+1, 6j] & \mbox{if } s=8j+2;\\
[2j+1, 6j+1] & \mbox{if }s=8j+3,\end{cases}
\Longleftrightarrow t\in \begin{cases} [k, 3k-1] & \mbox{if }s=4k;\\
[k, 3k] & \mbox{if } s=4k+1;\\
[k+1, 3k] & \mbox{if } s=4k+2;\\
[k+1, 3k+1] & \mbox{if }s=4k+3.\end{cases}\]
where $j,k\ge 1$.

\ms\nt We now show that $\th_{s} = \th(2t, 4s-6-2t,2s-4,(4s-6)^{[s-3]})$,   for $s\ge 4$ and $\frac{2s-3}{8}< t< \frac{6s-5}{8}$, admits a local antimagic 2-coloring.
We keep the notation defined above. Following is a general approach:
\begin{enumerate}[{Step }1:]
\item Label the edges of the path $R_j$ of length $4s-6$ by the sequence\\ $A_{2s-3}(j; 2s-3)\diamond A_{2s-3}(x-j; -(2s-3))$ in order for $1\le j\le s-2$.

\item For convenience, write $x_{s-2}=\a$. Separate $R_{s-2}$ into two paths. The first $2t$ edges form the path $Q_2$ and the rest form the path $Q_1$. So $\a$ and $\g_1$ are labeled at the end-edges incident to $u$. Recall that $\g_1=t(2s-3)+\a$.

\item Label the edges of the $(u,v)$-path $P$ of length $2s-4$ by the reverse of the sequence $A_{s-2}(2s-3; 2s-3)\diamond A_{s-2}((2s-3)(2s-4); -2s+3)$,  i.e., $A_{s-2}((s-1)(2s-3); 2s-3)\diamond A_{s-2}((s-2)(2s-3); -2s+3)$.
\end{enumerate}
\nt Clearly, by the construction above, it induces a local antimagic labeling for $\th(2t, 4s-6-2t,2s-4,(4s-6)^{[s-3]})$. Under this labeling, the induced vertex label for $u$ is
\[(s-1)(2s-3) +\sum\limits_{i=1}^{s-2}i +\gamma_1=(2s-3)(s-1+t)+\frac{s^2-3s+2}{2}+\a.\]
The difference from $y=(2s-3)(2s-2)$ is $\delta(t)=(2s-3)(s-1-t)-\frac{s^2-3s+2}{2}-\a$. Clearly $\delta(t)$ is a decreasing function of $t$.

\ms\nt Now, if we choose $\a=1$, then $\delta(t)=\frac{3s^2-7s-4st+6t+2}{2}$, where $\frac{2s-3}{8}< t<\frac{6s-5}{8}$.
So
\[\left.\begin{array}{r}
16k^2-11k+1\\
16k^2-k-1\\
16^2+k-1\\
16k^2+11k+1
\end{array}\right\}\ge\delta(t)\ge \begin{cases}
3k-2 & \mbox{if }s=4k;\\
k-1 & \mbox{if }s=4k+1;\\
7k & \mbox{if }s=4k+2;\\
5k+1 & \mbox{if }s=4k+3.
\end{cases}\]

\nt The set of differences of two end-edge labels in $R_j$, $2\le j \le s-2$, is  $D=\{1,3,\ldots,2s-7\}=A_{s-3}(1;2)$.

\nt Clearly $\delta(t)=2$ only when $(s,t)=(13,9)$. Also the maximum value of $\delta(t)$ for each case of $s$ is greater than $(s-3)^2$. Let us look at the second and third largest values $\delta_2$ and $\delta_3$ of $\delta(t)$ if any:
\[\delta_2=\begin{cases}
16k^2-19k+4 & \mbox{if }s=4k;\\
16k^2-9k & \mbox{if }s=4k+1;\\
16k^2-7k-2 & \mbox{if }s=4k+2;\\
16k^2+3k-2 & \mbox{if }s=4k+3.
\end{cases}\qquad \delta_3=\begin{cases}
16k^2-27k+7 & \mbox{if }s=4k;\\
16k^2-17k+1 & \mbox{if }s=4k+1;\\
16k^2-15k-3 & \mbox{if }s=4k+2;\\
16k^2-5k-5 & \mbox{if }s=4k+3.
\end{cases}\]
Clearly $0\le \delta_3<(s-3)^2-2$. So by Lemma~\ref{lem-AP-2}, there is a subset $B$ of $D$ such that the sum of integers in $B$ is $\delta(t)$ when $\frac{2s-3}{8}+2< t<\frac{6s-5}{8}$ except the cases $(s,t)=(13,9)$. Similar to Case (2), we find a local antimagic $2$-coloring for $\th(2t, 4s-6-2t,2s-4,(4s-6)^{[s-3]})$ according to the above range of $t$.

\nt For the case $(s,t)=(13, 9)$, $y=552$. Under the proposed labeling we can see that the induced label for $u$ is $549+\a$. So we may choose $\a=3$.

\ms\nt The remaining cases is when $\frac{2s-3}{8}< t\le\frac{2s-3}{8}+2$. When $s=4$, we have $\delta_2=1$ and $\delta_3$ does not exist.
 We shall modify our proposed labeling. Now, we choose $\a=2s-4$. In this case, $1$ is not labeled at the end-edge incident to $u$ so that the set of labels of the end-edges incident to $u$ is $\{(s-1)(2s-3),\gamma_1\}\cup [2, s-2]\cup\{2s-4\}$. Thus, the sum is $(s-1)(2s-3) +(2s-4)+\sum\limits_{i=2}^{s-2}i +\gamma_1=(2s-3)(s-1+t)+\frac{s^2+5s-16}{2}$. The difference from $y=(2s-3)(2s-2)$ is $\delta^*(t)=\frac{3s^2-15s-4st+6t+22}{2}$. One may easily check that $3\le \delta^*(t)\le (s-3)^2-3$ for $\frac{2s-3}{8}< t\le\frac{2s-3}{8}+2$, except $(s,t)=(4,2), (5,2), (6,3), (7,3)$. Thus we have a local antimagic $2$-coloring for $\th(2t, 4s-6-2t,2s-4,(4s-6)^{[s-3]})$ when $\frac{2s-3}{8}< t\le\frac{2s-3}{8}+2$.

\ms\nt For those exceptional cases, we have
\begin{enumerate}[1.]
\item $(s,t)=(4,2)$. Now $\delta(2)=1$. We may apply the original approach.
\item $(s,t)=(5,2)$. $\th_5=\th(4,6,10,14,14)$ with edge labels \\ $39,10,46,3;$ \\ $7,42,14,35,21,28;$ \\ $4,45,11,38,18,31,25,24,32,17;$ \\ $1,48,8,41,15,34,22,27,29,20,36,13,43,6;$ \\ $5,44,12,37,19,30,26,23,33,16,40,9,47,2$.
\item $(s,t)=(6,3)$. Now $\delta(3)=7<3^2$. We may apply the original approach.
\item $(s,t)=(7,3)$. Now $x=121$, $y=132$. $\th(6,10,16,22,22,22,22)$ with sequences \\
$4,117,15,106,26,95;$ \\
$66, 55, 77, 44, 88, 33, 99, 22, 110, 11;$\\ $37,84,48,72,59,62,70,51,81,40,92,29,103,18,114,7;$\\
$2,119,13,108,24,97,35,86,46,75,57,64,68,53,79,42,90,31,101,20,112,9;$\\
$5,116,16,105,27,94,38,83,49,72,60,61,71,50,82,39,93,28,104,17,115,6;$\\
$8, 113, 19, 102, 30, 91, 41, 80, 52, 69, 63, 58, 74, 47, 85, 36, 96, 25, 107, 14, 118, 3;$\\
$10, 111, 21, 100, 32, 89, 43, 78, 54, 67, 65, 56, 76, 45, 87, 34, 98, 23, 109, 12, 120, 1.$\\
\end{enumerate}

\nt So we have a local antimagic $2$-coloring for $\th(2t, 4s-6-2t,2s-4,(4s-6)^{[s-3]})$   when $s\ge 4$ and $\frac{2s-3}{8}< t< \frac{6s-5}{8}$.

\ms\nt Note that, one may see from each case that $m>2s+2$. This completes the proof.
\end{proof}

\section{Examples} In this section, we shall provide example(s) to illustrate the construction of each case and also provide solutions for the exceptional cases raised in the proof of Theorem~\ref{thm-chilath=2}.

\begin{example} The aim of this example is to illustrate the construction showed in Case (1).

\ms\nt Take $s=6$ (i.e., $k=1$), we have $\th_6=\th(4,4,4,4,4,6)$ with $m=26$, $x=27$, $y=39$, $U_1=\{1\}$, $U_2 = \{4,5,8,9,12\}$, $[1,12]\setminus(U_1\cup U_2)=\{2,3,6,7,10,11\}$.

\nt $A_3(1; 12)=(1, 13, 25)$ and $A_3(26; -12)=(26, 14, 2)$.   So $A_3(1; 12) \diamond A_3(26; -12)=(1, 26, 13, 14, 25, 2)$.

\nt Similarly,

\nt $A_2(4; 12)=(4, 16)$ and $A_2(23, -12)=(23, 11)$, $A_2(5; 12)=(5,17)$ and $A_2(22; -12)=(22, 10)$, $A_2(8; 12)=(8, 20)$ and $A_2(19;-12)=(19, 7)$, $A_2(9; 12)=(9, 21)$ and $A_2(18; -12)=(18, 6)$, $A_2(12; 12)=(12, 24)$ and $A_2(15; -12)=(15,3)$.

\nt So, the paths of length 4 and 6 have edge labels \[4,23,16,11; \quad  5,22,17,10; \quad 8,19,20,7; \quad 9,18,21,6; \quad 12,15,24,3; \quad 1,26,13,14,25,2.\]  All the left (respectively right) end vertices are merged to get the degree 6 vertex with induced label 39.  \rsq
\end{example}

\begin{example} The aim of this example is to illustrate the construction showed in Case~(2).

\ms\nt Take $s=9$ (i.e., $l=3$), we get $\th(4,10^{[8]})$ with $y=102$, $x=85$. Keep the notation defined in Lemma~\ref{lem-AP-2} and the proof of Theorem~\ref{thm-chilath=2}. Since $\delta = 15$, $n=8$, we choose $\kappa= 15$ with $\tau=0$. By Lemma~\ref{lem-AP-2}, we have $B=\{15\}$. So we replace 1 by 16 as a label of end-edge incident to $u$. Thus $u$ is incident to end-edge labels in $\{16, 2, 3, 4, 5, 6, 7, 8, 51\}$. The paths labels are \\
$51,34,68,17$: $A_2(51;17)\diamond A_2(34;-17)$;\\
$16,69,33,52,50,35,67,18,84,1$: the {\it reverse of} $A_5(1;17)\diamond A_5(84;-17)$;\\
$2,83,19,66,36,49,53,32,70,15$: $A_5(2;17)\diamond A_5(83;-17)$;\\
$3,82,20,65,37,48,54,31,71,14$: $A_5(3;17)\diamond A_5(82;-17)$;\\
$4,81,21,64,38,47,55,30,72,13$: $A_5(4;17)\diamond A_5(81;-17)$;\\
$5,80,22,63,39,46,56,29,73,12$: $A_5(5;17)\diamond A_5(80;-17)$;\\
$6,79,23,62,40,45,57,28,74,11$: $A_5(6;17)\diamond A_5(79;-17)$;\\
$7,78,24,61,41,44,58,27,75,10$: $A_5(7;17)\diamond A_5(78;-17)$;\\
$8,77,25,60,42,43,59,26,76,9$: $A_5(8;17)\diamond A_5(77;-17)$.\\

\ms\nt Using $s=12$ (i.e., $l=4$), we get $\th(6,14^{[11]})$ with $y=184$, $x=161$.   Since $\delta = 26$. We choose $\kappa=21$ (i.e., $k=1$) with $\tau=5$. By Lemma~\ref{lem-AP-2} we have $B=\{21, 5\}$. So we replace 1 by 22 and 9 by 14 as labels of end-edges incident to $u$. Thus $u$ is incident to end-edge labels in $\{22, 2, 3, 4, 5, 6, 7, 8, 14, 10, 11 ,92\}$. The paths labels are \\
$92, 69, 115, 46, 138, 23$: $A_3(92;23)\diamond A_3(69;-23)$;\\
$22,139,45,116,68,93,91,70,114,47,137,24,160,1$: the reverse of $A_7(1;23)\diamond A_7(160;-23)$;\\
$2,159,25,136,48,113,71,90,94,67,117,44,140,21$: $A_7(2;23)\diamond A_7(159;-23)$;\\
$3,158,26,135,49,112,72,89,95,66,118,43,141,20$: $A_7(3;23)\diamond A_7(158;-23)$; \\
$4,157,27,134,50,111,73,88,96,65,119,42,142,19$: $A_7(4;23)\diamond A_7(157;-23)$; \\
$5,156,28,133,51,110,74,87,97,64,120,41,143,18$: $A_7(5;23)\diamond A_7(156;-23)$; \\
$6,155,29,132,52,109,75,86,98,63,121,40,144,17$: $A_7(6;23)\diamond A_7(155;-23)$; \\
$7,154,30,131,53,108,76,85,99,62,122,39,145,16$: $A_7(7;23)\diamond A_7(154;-23)$; \\
$8,153,31,130,54,107,77,84,100,61,123,38,146,15$: $A_7(8;23)\diamond A_7(153;-23)$; \\
$14,147,37,124,60,101,83,78,106,55,129,32,152,9$: the reverse of $A_7(9;23)\diamond A_7(152;-23)$; \\
$10,151,33,128,56,105,79,82,102,59,125,37,148,13$: $A_7(10;23)\diamond A_7(151;-23)$; \\
$11,150,34,127,57,104,80,81,103,58,126,36,149,12$: $A_7(11;23)\diamond A_7(150;-23)$.
\rsq
\end{example}

\begin{example}\label{ex-3-1} The aim of this example is to illustrate the construction showed in Case~(3) and provide a local antimagic $2$-coloring for the exceptional case $(l,t)=(6,7)$.

\nt
Let $s=12$, i.e., $l=6$. Now, $x=231$ and $y=253$.

\begin{enumerate}[(a)]

\item The graph is $\th_{12}=\th(22-2t, 2t, 20^{[6]}, 22^{[4]})$, where $t=6,7$.
Begin with the sequences\\
{\fontsize{8}{9}\selectfont
$A_{11}(6; 22)\diamond A_{11}(225; -22)$: 6, 225, 28, 203, 50, 181, 72, 159, 94, 137, 116, 115, 138, 93, 160, 71, 182, 49, 204, 27, 226, 5\\
$A_{11}(7; 22)\diamond A_{11}(224; -22)$: 7, 224, 29, 202, 51, 180, 73, 158, 95, 136, 117, 114, 139, 92, 161, 70, 183, 48, 205, 26, 227, 4\\
$A_{11}(8; 22)\diamond A_{11}(223; -22)$: 8, 223, 30, 201, 52, 179, 74, 157, 96, 135, 118, 113, 140, 91, 162, 69, 184, 47, 206, 25, 228, 3\\
$A_{11}(9; 22)\diamond A_{11}(222; -22)$: 9, 222, 31, 200, 53, 178, 75, 156, 97, 134, 119, 112, 141, 90, 163, 68, 185, 46, 207, 24, 229, 2\\
$A_{11}(10; 22)\diamond A_{11}(221; -22)$: 10, 221, 32, 199, 54, 177, 76, 155, 98, 133, 120, 111, 142, 89, 164, 67, 186, 45, 208, 23, 230, 1

$A_{10}(11; 22)\diamond A_{10}(220; -22)$: 11, 220, 33, 198, 55, 176, 77, 154, 99, 132, 121, 110, 143, 88, 165, 66, 187, 44, 209, 22\\
$A_{10}(12; 22)\diamond A_{10}(219; -22)$: 12, 219, 34, 197, 56, 175, 78, 153, 100, 131, 122, 109, 144, 87, 166, 65, 188, 43, 210, 21\\
$A_{10}(13; 22)\diamond A_{10}(218; -22)$: 13, 218, 35, 196, 57, 174, 79, 152, 101, 130, 123, 108, 145, 86, 167, 64, 189, 42, 211, 20\\
$A_{10}(14; 22)\diamond A_{10}(217; -22)$: 14, 217, 36, 195, 58, 173, 80, 151, 102, 129, 124, 107, 146, 85, 168, 63, 190, 41, 212, 19\\
$A_{10}(15; 22)\diamond A_{10}(216; -22)$: 15, 216, 37, 194, 59, 172, 81, 150, 103, 128, 125, 106, 147, 84, 169, 62, 191, 40, 213, 18\\
$A_{10}(16; 22)\diamond A_{10}(215; -22)$: 16, 215, 38, 193, 60, 171, 82, 149, 104, 127, 126, 105, 148, 83, 170, 61, 192, 39, 214, 17
}

Now the difference sets are $D_1=A_5(-1;-2)$ and $D_2=A_6(1;2)$.
\begin{enumerate}[i)]
\item $t=6$. So $\th_{12}=\th(10, 12, 20^{[6]}, 22^{[4]})$. Initially, we use the first five sequences above to label the $(u,v)$-paths $T_j$ and the last six sequences above to label the $(u,v)$-paths $R_i$. We then break $T_5$ into two parts such that the first 10 edges form the $(u,v)$-path $Q_2$ and the remaining 12 edges form the $(u,v)$-path $Q_1$. Now, the induced vertex label for $u$ is $\sum\limits_{j=6}^{16} j+ 120=241$. Thus $\d(6)=12$. So we choose $B=\{1,11\}\subset D_2$. Therefore, the actual assignment for each $(u,v)$-path is to label:\\
    $T_1$ by $A_{11}(6; 22)\diamond A_{11}(225; -22)$; $T_2$ by $A_{11}(7; 22)\diamond A_{11}(224; -22)$; $T_3$ by $A_{11}(8; 22)\diamond A_{11}(223; -22)$; $T_4$ by $A_{11}(9; 22)\diamond A_{11}(222; -22)$;\\
    $Q_2$ by 10, 221, 32, 199, 54, 177, 76, 155, 98, 133;\\
    $Q_1$ by 120, 111, 142, 89, 164, 67, 186, 45, 208, 23, 230, 1;\\
    $R_1$ by the reverse of $A_{10}(11; 22)\diamond A_{10}(220; -22)$; $R_2$ by $A_{10}(12; 22)\diamond A_{10}(219; -22)$; $R_3$ by $A_{10}(13; 22)\diamond A_{10}(218; -22)$; $R_4$ by $A_{10}(14; 22)\diamond A_{10}(217; -22)$; $R_5$ by $A_{10}(15; 22)\diamond A_{10}(216; -22)$; $R_6$ by the reverse of $A_{10}(16; 22)\diamond A_{10}(215; -22)$. \\
    Thus,
    \[f^+(u)=6+7+8+9+10+120+22+12+13+14+15+17=253.\]

\item $t=7$. So $\th_{12}=\th(8, 14, 20^{[6]}, 22^{[4]})$. Initially, we use the first five sequences above to label the $(u,v)$-paths $T_j$ and the last six sequences above to label the $(u,v)$-paths $R_i$. We then break $T_5$ into two parts such that the first 8 edges form the $(u,v)$-path $Q_2$ and the remaining 14 edges form the $(u,v)$-path $Q_1$. Now, the induced vertex label for $u$ is $\sum\limits_{j=6}^{16} j+ 98=219$. Thus $\d(7)=34$.
   For this case, we do not have $B\subset D_2$. So we choose $B=\{-1,3,5,7,9,11\}\subset D_1\cup D_2$. Thus the actual assignment for each $(u,v)$-path is to label:\\
    $T_1$ by the reverse of $A_{11}(6; 22)\diamond A_{11}(225; -22)$; $T_2$ by $A_{11}(7; 22)\diamond A_{11}(224; -22)$; $T_3$ by $A_{11}(8; 22)\diamond A_{11}(223; -22)$; $T_4$ by $A_{11}(9; 22)\diamond A_{11}(222; -22)$;\\
    $Q_2$ by 10, 221, 32, 199, 54, 177, 76, 155;\\
    $Q_1$ by 98, 133, 120, 111, 142, 89, 164, 67, 186, 45, 208, 23, 230, 1;\\
    $R_1$ by the reverse of $A_{10}(11; 22)\diamond A_{10}(220; -22)$; $R_2$ by the reverse of $A_{10}(12; 22)\diamond A_{10}(219; -22)$; $R_3$ by the  reverse of $A_{10}(13; 22)\diamond A_{10}(218; -22)$; $R_4$ by the reverse of $A_{10}(14; 22)\diamond A_{10}(217; -22)$; $R_5$ by the reverse of $A_{10}(15; 22)\diamond A_{10}(216; -22)$; $R_6$ by $A_{10}(16; 22)\diamond A_{10}(215; -22)$. \\Thus,
    \[f^+(u)=5+7+8+9+10+98+22+21+20+19+18+16=253.\]
\end{enumerate}
\item The graph is $\th_{12}=\th(22-2t, 2t-2, 20^{[5]}, 22^{[5]})$, where $t=6,7$.
We begin with the following sequences that are the reverse of the initial sequences in Case (a): $A_{11}(1; 22)\diamond A_{11}(230;-22)$, $A_{11}(2; 22)\diamond A_{11}(229;-22)$, $A_{11}(3; 22)\diamond A_{11}(228;-22)$, $A_{11}(4; 22)\diamond A_{11}(227;-22)$, $A_{11}(5; 22)\diamond A_{11}(226;-22)$, $A_{10}(17; 22)\diamond A_{10}(214;-22)$, $A_{10}(18; 22)\diamond A_{10}(213;-22)$, $A_{10}(19; 22)\diamond A_{10}(212;-22)$, $A_{10}(20; 22)\diamond A_{10}(211;-22)$, $A_{10}(21; 22)\diamond A_{10}(210;-22)$, $A_{10}(22; 22)\diamond A_{10}(209;-22)$.

 Now, the difference sets are $D_1=A_5(1;2)$ and $D_2=A_6(-1,-2)$.
 \begin{enumerate}[i)]
\item $t=6$. So $\th_{12}=\th(10, 10, 20^{[5]}, 22^{[5]})$. Initially, we use the first five sequences above to label the $(u,v)$-paths $T_j$ and the last six sequences above to label the $(u,v)$-paths $R_i$.  We then break $R_6$ into two parts such that the first 10 edges form the $(u,v)$-path $Q_2$ and the remaining 10 edges form the $(u,v)$-path $Q_1$. Now, the induced vertex label of $u$ is $\sum\limits_{j=1}^{5} j+\sum\limits_{i=17}^{22} i +132 =264$. So we choose $B=\{-9,-3,1\}\subset D_1\cup D_2$.

    Thus the actual assignment for each $(u,v)$-path is to label:\\
    $T_1$ by $A_{11}(1; 22)\diamond A_{11}(230; -22)$; $T_2$ by $A_{11}(2; 22)\diamond A_{11}(229; -22)$; $T_3$ by $A_{11}(3; 22)\diamond A_{11}(228; -22)$; $T_4$ by $A_{11}(4; 22)\diamond A_{11}(227; -22)$; $T_5$ by the reverse of $A_{11}(5; 22)\diamond A_{11}(226; -22)$;\\
    $R_1$ by $A_{10}(17; 22)\diamond A_{10}(214; -22)$; $R_2$ by the reverse of $A_{10}(18; 22)\diamond A_{10}(213; -22)$; $R_3$ by $A_{10}(19; 22)\diamond A_{10}(212; -22)$; $R_4$ by $A_{10}(20; 22)\diamond A_{10}(211; -22)$; $R_5$ by the reverse of $A_{10}(21; 22)\diamond A_{10}(210; -22)$;\\
    $Q_2$ by 22, 209, 44, 187, 66, 165, 88, 143, 110, 121;\\
    $Q_1$ by 132, 99, 154, 77, 176, 55, 198, 33, 220, 11. \\
    Thus,
    \[f^+(u)=1+2+3+4+6+17+15+19+20+12+22+132=253.\]
\item $t=7$. So $\th_{12}=\th(8, 12, 20^{[5]}, 22^{[5]})$.
Initially, we use the first five sequences above to label the $(u,v)$-paths $T_j$ and the last six sequences above to label the $(u,v)$-paths $R_i$.  We then break $R_6$ into two parts such that the first 8 edges form the $(u,v)$-path $Q_2$ and the remaining 12 edges form the $(u,v)$-path $Q_1$. Now, the induced vertex label of $u$ is $\sum\limits_{j=1}^{5} j+\sum\limits_{i=17}^{22} i +110 =242$. Now $\d(6)=11$. So we may choose $B=\{1,3,7\}$.

Thus the actual assignment for each $(u,v)$-path is to label:\\
    $T_1$ by $A_{11}(1; 22)\diamond A_{11}(230; -22)$; $T_2$ by the reverse of $A_{11}(2; 22)\diamond A_{11}(229; -22)$; $T_3$ by $A_{11}(3; 22)\diamond A_{11}(228; -22)$; $T_4$ by the reverse of $A_{11}(4; 22)\diamond A_{11}(227; -22)$; $T_5$ by the reverse of $A_{11}(5; 22)\diamond A_{11}(226; -22)$;\\
    $R_1$ by $A_{10}(17; 22)\diamond A_{10}(214; -22)$; $R_2$ by $A_{10}(18; 22)\diamond A_{10}(213; -22)$; $R_3$ by $A_{10}(19; 22)\diamond A_{10}(212; -22)$; $R_4$ by $A_{10}(20; 22)\diamond A_{10}(211; -22)$; $R_5$ by $A_{10}(21; 22)\diamond A_{10}(210; -22)$;\\
    $Q_2$ by 22, 209, 44, 187, 66, 165, 88, 143;\\
    $Q_1$ by 110, 121, 132, 99, 154, 77, 176, 55, 198, 33, 220, 11. \\
    Thus,
    \[f^+(u)=1+9+3+7+6+17+18+19+20+21+22+110=253.\]
 \end{enumerate} \rsq
\end{enumerate}
\end{example}

\begin{example}\label{ex-3-2} The aim of this example is to illustrate the construction showed in Case~(3) and provide a local antimagic $2$-coloring for the exceptional case $(l,t)=(3,3)$.
Let $s=6$, i.e., $l=3$. Now, $x=45$ and $y=55$. The sequences are\\
$A_5(1; 10)\diamond A_5(44; -10)$: 1, 44, 11, 34, 21, 24, 31, 14, 41, 4\\
$A_5(2; 10)\diamond A_5(43; -10)$: 2, 43, 12, 33, 22, 23, 32, 13, 42, 3\\
$A_4(5; 10)\diamond A_4(40; -10)$: 5, 40, 15, 30, 25, 20, 35, 10\\
$A_4(6; 10)\diamond A_4(39; -10)$: 6, 39, 16, 29, 26, 19, 36, 9\\
$A_4(7; 10)\diamond A_4(38; -10)$: 7, 38, 17, 28, 27, 18, 37, 8

\begin{enumerate}[(a)]
\item $t=l=3$. So $\th_6=\th(4, 6, 8^{[3]}, 10)$.

$(u,v)$-path $T_1$ is labeled by 4, 41, 14, 31; 24, 21, 34, 11, 44, 1. So\\
$(u,v)$-path $Q_2$ is labeled by 4, 41, 14, 31 and\\ $(u,v)$-path $Q_1$ is labeled by 24, 21, 34, 11, 44, 1.

$(u,v)$-path $T_2$ is labeled by 3, 42, 13, 32, 23, 22, 33, 12, 43, 2.

$(u,v)$-path $R_1$ is labeled by 10, 35, 20, 25, 30, 15, 40, 5.

$(u,v)$-path $R_3$ is labeled by 8, 37, 18, 27, 28, 17, 38, 7.

$(u,v)$-path $R_2$ is labeled by 6, 39, 16, 29, 26, 19, 36, 9

Thus, $f^+(u)=4+24+3+10+ 8+6 =55$.

\item $t=l=3$. So $\th_6=\th(4, 4, 8^{[2]},10^{[2]})$.

$(u,v)$-path $Q_2$ is labeled by 8, 37, 18, 27.

$(u,v)$-path $Q_1$ is labeled by 28, 17, 38, 7.

$(u,v)$-path $R_1$ is labeled by 6, 39, 16, 29, 26, 19, 36, 9.

$(u,v)$-path $R_2$ is labeled by 10, 35, 20, 25, 30, 15, 40, 5.

$(u,v)$-path $T_1$ is labeled by 1, 44, 11, 34, 21, 24, 31, 14, 41, 4.

$(u,v)$-path $T_2$ is labeled by 2, 43, 12, 33, 22, 23, 32, 13, 42, 3.

Thus, $f^+(u)=8+28+6+10+ 1+2 =55$. \rsq
\end{enumerate}

\end{example}

\begin{example} The aim of this example is to illustrate the construction given in Case~(4).
 Take $s=7$ so that $\th_7=\th(2t, 22-2t, 10, 22^{[4]})$, $2\le t\le 4$.  We have $x=121$, $y=132$ and $y-x=11$.

{\fontsize{9}{9}\selectfont
$\begin{aligned}
A_{11}(1;11)\diamond A_{11}(120;-11)& = 1, 120, 12, 109, 23, 98, 34, 87, 45, 76, 56, 65, 67, 54, 78, 43, 89, 32, 100, 21, 111, 10; &\\
A_{11}(2;11)\diamond A_{11}(119;-11) & = 2, 119, 13, 108, 24, 97, 35, 86, 46, 75, 57, 64, 68, 53, 79, 42, 90, 31, 101, 20, 112, 9;& [7]\\
A_{11}(3;11)\diamond A_{11}(118;-11) &  = 3, 118, 14, 107, 25, 96, 36, 85, 47, 74, 58, 63, 69, 52, 80, 41, 91, 30, 102, 19, 113, 8; & [5]\\
A_{11}(4;11)\diamond A_{11}(117;-11) & = 4, 117, 15, 106, 26, 95, 37, 84, 48, 72, 59, 62, 70, 51, 81, 40, 92, 29, 103, 18, 114, 7; & [3]\\
A_{11}(5;11)\diamond A_{11}(116;-11) & = 5, 116, 16, 105, 27, 94, 38, 83, 49, 72, 60, 61, 71, 50, 82, 39, 93, 28, 104, 17, 115, 6. &[1]\\
A_{5}(66;11)\diamond A_5(55;-11) & = 66, 55, 77, 44, 88, 33, 99, 22, 110, 11 \leftarrow \mbox{ this sequence is for the $(u,v)$-path $P$.}
\end{aligned}$}

Note that $(s-3)^2=16$. The number with a bracket behind the sequence is the difference between the last and the first terms. Hence $D=\{1,3,5,7\}$.
\begin{enumerate}[1.]
\item When $t=4$. We have $\delta(4)=6<16$. First we separate $A_{11}(1;11)\diamond A_{11}(120;-11)$ into two sequences: 1, 120, 12, 109, 23, 98, 34, 87; and  45, 76, 56, 65, 67, 54, 78, 43, 89, 32, 100, 21, 111, 10.
Since $\delta(4)<7$, by Lemma~\ref{lem-AP-2}, we choose $B=\{1,5\}$. So we reverse the order of $A_{11}(5;11)\diamond A_{11}(116;-11)$ and $A_{11}(3;11)\diamond A_{11}(118;-11)$, i.e., the end-edge labels for $u$ is 1, $45=\gamma_1$, 2, 8, 4, 6, 66.
\item When $t=3$. We have $\delta(3)=17>16$ and $\delta^*(3)=-1$. We must use an ad hoc method which is shown in the proof.
\item When $t=2$. We have $\delta(2)=28>16$. $\delta^*(2)=10<16$. First we separate the reverse of $A_{11}(1;11)\diamond A_{11}(120;-11)$ into two sequences: 10, 111, 21, 100; and 32, 89, 43, 78, 54, 67, 65, 56, 76 45, 87, 34, 98, 23, 109, 12, 120, 1.
Since $\delta^*(2)=10$, we choose $B=\{7, 3\}$. So we reverse the order of $A_{11}(2;11)\diamond A_{11}(119;-11)$ and $A_{11}(4;11)\diamond A_{11}(117;-11)$, i.e., the end-edge labels for $u$ is 10, $32=\gamma_1$, 9, 3, 7, 5, 66. \rsq
\end{enumerate}
\end{example}

\section{Conjecture and Open Problem}

We have completely characterized $s$-bridge graphs $\th_s$ with $\chi_{la}(\th_s)=2$. We note that the only other known results on $s$-bridge graphs are (i) $\chi_{la}(\th(a,b))=3$ for $a,b\ge 1$ and $a+b\ge 3$; and (ii) $\th(2^{[s]})=3$ for odd $s\ge 3$. We end with the following conjecture and open problem.

\begin{conjecture} If $\th_s$ is not a graph in Theorem~\ref{thm-chilath=2}, then $\chi_{la}(\th_s) = 3$. \end{conjecture}

\begin{problem} Characterize graph $G$ with $\chi_{la}(G)=2$. \end{problem}






\end{document}